\documentclass[11pt,a4paper]{article} 
\usepackage[cp866]{inputenc}
\usepackage[english]{babel}
\usepackage{amsmath}
\usepackage{amsfonts}
\usepackage{longtable}
\usepackage{amssymb}
\usepackage{graphics}
\newcommand{\nid}{\noindent}
\newcommand{\bc}{\begin{center}}
\newcommand{\br}{\begin{right}}
\newcommand{\ec}{\end{center}}
\newcommand{\be}{\begin{equation}}
\newcommand{\ee}{\end{equation}}

\newcommand{\vl}{\mid}

\newcommand{\rar}{\rightarrow}
\newcommand{\grad}{\nabla}

\newcommand{\p}{\partial}

\newcommand{\bint}{\mbox{int} \,}

\newcommand{\bco}{\mbox{co} \,}

\newcommand{\exs}{\exists}

\newcommand{\meq}{\geq}

\newcommand{\al}{\alpha}

\newcommand{\eps}{\varepsilon}

\newtheorem{thm}{Theorem}[section]

\newtheorem{rem}{Remark}[section]
\newtheorem{defi}{Definition}[section]
\newtheorem{cor}{Corollary}[section]
\newtheorem{ex}{Example}[section]

\newtheorem{lem}{Lemma}[section]
\newcommand{\avin}[6]{#3^{-1} \; \int^#2_#1 \, \grad #6(r(#4, \tau, #5))d\tau}

\newcommand{\incl}{\subset}

\begin{document}

\hspace{5.8cm} Devoted to my teacher Prof. V.F. Demyanov

\hspace{5.8cm} who formulated the problem about the sub-

\hspace{5.8cm} differential of the second order long time ago

\vspace{0.5cm}

AMS 517.9

\vspace{0.5cm}

\bc {\bf Prudnikov I.M.} \ec \vspace{0.5cm}
\begin{center}
{\large \bf THE SUBDIFFERENTIALS OF THE FIRST AND SECOND ORDERS
FOR LIPSCHITZ FUNCTIONS}
\end{center}

\vspace{0.5cm}

Construction of an united theory of the subdifferentials of the
first and second orders is interesting for many specialists in
optimization \cite{morduchrockafel}. In the paper the rules for
construction  of the  subdifferentials of the first and second
orders are introduced. The constructions are done with the help of
the Steklov integral of a Lipschitz function $f(\cdot)$ over the
images of a set-valued mapping $D(\cdot)$. It is proved that the
subdifferential of the first order  consisting of the average
integral limit values of the gradients $\nabla f(r(\cdot))$,
calculated along the curves $r(\cdot)$ from an introduced set of
curves $\eta$, coincides with the subdifferentials of the first
order constructed using the Steklov integral introduced by the
author for the first time in \cite{lowapp2}, \cite{lowapp2a}. If
the function $f(\cdot)$ is twice differentiable at $x$ then  the
subdifferentials of the first and second orders coincide with the
gradient $\nabla f(r(\cdot))$ and the matrix of the second mixed
derivatives of $f(\cdot)$ at $x$. The generalized gradients and
matrices are used for formulation of the necessary and sufficient
conditions of optimality. The calculus for the  subdifferentials
of the first and second orders is constructed. The examples are
given.

\nid


\bigskip
\normalsize \section{Introduction}

Lipschitz functions are not smooth in general case. They are
almost everywhere (a.e.) differentiable in $\mathbb{R}^n$. The
generalized gradients are introduced for formulation of the
necessary conditions of optimality. Union over them is called the
subdifferential. There is not the unique way for introduction of
the subdifferential. So F. Clarke introduced the subdifferential
consisting of limit values of the gradients of  Lipschitz
function. This subdifferential was called the Clarke
subdifferential \cite{clark1}, \cite{clark}. Mischel and Penot
defined their subdifferential and generalized gradients by means
of the limit ratio of values of function at points from a
neighborhood of a considered point and an increment of argument
\cite{penot2}. The necessary optimality condition in
$\mathbb{R}^n$ can be written in the form that zero belongs to the
subdifferential. In smooth case this condition can be rewritten in
the form that the derivative is equal to zero.

The author introduced the new subdifferential that always belongs
to the Clarke subdifferential. If a function is difference of two
convex functions (so called DC functions) then this
subdifferential is equal to the Clarke subdifferential
\cite{lowapp2}. The subdifferential consists of the average
integral limit values of the gradients of function, calculated
along curves from a set of curves, introduced in \cite{lowapp2},
\cite{lowapp2a}. If function is differentiable at $x$ then the
defined subdifferential is equal to the gradient of $f(\cdot)$ at
$x$ i.e. $\nabla f(x)= f'(x)$. The necessary optimality condition
at $x$ can be written in the form: zero belongs to the
subdifferential calculated at $x$.

To write a condition of the second order of optimality it is
necessarily to introduce the subdifferential of the second order,
consisting of the generalized matrices. Many mathematicians tried
to introduce the subdifferential of the second order in different
ways. So, in \cite{morduchrockafel} the authors introduced the
partial second-order subdifferential. Instead of this the
generalized matrices are considered in the paper, that are used
for formulation of the second order conditions of optimality.

 The difficulty for introduction the subdifferential of
the second order is that any Lipschitz function is not a.e. twice
differentiable in general case. Consequently, the previous method,
used for introduction  of the subdifferential of the first order,
is not appropriate. We need to develop a new method for the
definition of the generalized matrices of the second mixed
derivatives that is not connected with differentiable qualities of
any Lipschitz function. This problem was solved in this paper.

\section{\bf  The construction of the subdifferential of the first order}

Let  $f(\cdot):\mathbb{R}^n \rar \mathbb{R}$ be a Lipschitz
function with the Lipschitz constant $L$. Our goal is to
investigate the differential qualities of the first and second
orders of the function $f(\cdot)$ with the help of the Steklov
integral.

Define the function $\varphi(\cdot):\mathbb{R}^n \rar \mathbb{R}$
\be \varphi(x)= \frac{1}{\mu(D(x))} \int_{D(x)} f(x+y)d y,
\label{firstsecsub} \ee where $D(\cdot): \mathbb{R}^n \rar
2^{\mathbb{R}^n}$  is a continuous set-valued mapping (SVM) in the
Hausdorff metric with convex compact images, $\mu(D(x)$ is the
measure of $D(x)$. The function $\varphi(\cdot)$ depends on the
chosen SVM $D(\cdot).$

The integrals ($\ref{firstsecsub}$) for the constant SVM $D(\cdot)
\equiv D$ is called the Steklov integrals. Its qualities were
studied in \cite{proudintegapp1}-\cite{krotovprohorovich}. It was
proved that $\varphi(\cdot) $ is the Lipschitz continuously
differentiable function with the Lipschitz derivative
$\varphi'(\cdot) $. We define the  Lipschitz constant of
$\varphi(\cdot)$ by $L(D)$.

Let us consider SVM $D(\cdot)$ satisfying the following
conditions.
\begin{enumerate}
\item $x_0 \in \bint (x + D(x))$ for all $x \in S, \,\,\, S
\subset \mathbb{R}^n,$  is a neighborhood of $x_0$;

\item the diameter of $D(x)$ which we denote by $diam \, D(x) =
d(D(x)),$ tends to zero as $x \rar x_0$ and satisfies the
inequality $d(D(x)) \leq k \| x - x_0 \|$ for some constant $k$;

\item for some sequence $\{ \eps_i \}, \eps_i \rar +0,$ as $i \rar
\infty$ SVM $D(\cdot)$ is constant for all $x$ from the set
$\eps_{2i+1} < \parallel x - x_0 \parallel < \eps_{2i}$;

\item the boundary of $D(x)$ for all $x \in S,\, x\neq x_0, $ is
defined by continuously differentiable function of $x$.

\end{enumerate}

We will consider SVM $D(\cdot)$ satisfying the above mentioned
conditions for any sequences $\{ \eps_i \}, \, \eps_i \rar +0,$
and constants $k$. Denote the defined set of SVM by $\Xi$.
$D(\cdot)$ is constant SVM, i.e. $D(x) \equiv D_{2i}$, for any
$x$, satisfying $\eps_{2i+1} < \parallel x - x_0 \parallel <
\eps_{2i}$. The derivative $\varphi'(\cdot)$ is the Lipschitz
function with a constant $L_{2i}(D_{2i})$ (see.
\cite{proudintegapp1}, \cite{proudintegapp}).

Define  for SVM $D(\cdot)$ the set
$$
\p \varphi_D(x_0)=\bco \{ v \in \mathbb{R}^n \vl v= \lim_{x_i \rar
x_0} \varphi'(x_i) \},
$$
where the points $x_i$ are taken from the regions of constancy of
SVM $D(\cdot)$.  $\p \varphi_D(x_0)$ is a convex compact set in
$\mathbb{R}^n.$ The boundedness of $\p \varphi_D(x_0)$ follows
from the inequalities written below. We have
$$
 \varphi'(x)= \frac{1}{\mu(D(x))} \int_{D(x)} f'(x+y)d y
$$
for any $x$ from the regions of constancy of SVM $D(\cdot)$
(\cite{proudintegapp1}, \cite{proudintegapp}).

Since $\| f'(x+y) \| \leq L $, then  the inequality
$$
\| \varphi'(x) \| \leq L.
$$
follows from here. Consequently, $\p \varphi_D(x_0)$ is a convex
bounded set. Let us prove its closure.

Take $v_i \in \p \varphi_{D}(x_0)$ and $v_i \rar v$ for SVM
$D(\cdot) \in \Xi$. We will prove, that $v \in \p
\varphi_{D}(x_0)$. Consider a  subsequence $\{ i_k \} \subset \{ i
\} $ such that for $j \in \{ i_k \}$ and for the points $x_j$ from
the regions of constancy of SVM $D(\cdot)$ the inequality
$$
\|  \frac{1}{\mu(D(x_j))} \int_{D(x_j)} f'(x_j+y)d y - v_j \| \leq
\eps_j,
$$
was true where $\eps_j \xrightarrow[j]{} +0. $ Going to the limit
in $j$ we will have
$$
\lim_{x_j \rar x_0} \varphi'(x_j) = \lim_{x_j \rar x_0}
\frac{1}{\mu(D(x_j))} \int_{D(x_j)} f'(x_j+y)d y   = v,
$$
i.e. $v \in \p \varphi_{D}(x_0)$. Consequently, $\p
\varphi_{D}(x_0)$ is a closed set. Compactness of
$\varphi_{D}(x_0)$ follows from its boundedness  and closure.

We have proved the lemma
\begin{lem}
$\p \varphi_{D} (x_0) $ is a convex compact set in $\mathbb{R}^n$.
\label{lemfirstsecsub0}
\end{lem}

We will prove further that $\varphi'(\cdot)$ is continuous at
points $x, x \neq x_0,$ for any SVM $D(\cdot)$ satisfying the
conditions written above.

Consider the case when $f(\cdot)$ is differentiable at $x_0$. Find
the vectors which the set $\p \varphi_D(x_0)$ consists of. We have
$$ \varphi(x_0+\triangle x)=
f(x_0)+(f'(x_0),\triangle x)+\frac{1}{\mu(D(x_0+\triangle x))}
\int_{D(x_0+\triangle x)} (f'(x_0),y) d y +
$$
\be + \frac{1}{\mu(D(x_0+\triangle x))} \int_{D(x_0+\triangle x)}
o(\triangle x + y) d y. \label{firstsecsub2} \ee Show that
$$
\frac{1}{\mu(D(x_0+\triangle x))} \int_{D(x_0+\triangle x)}
o(\triangle x + y) d y = \tilde{o}(\triangle x),
$$
where $\tilde{o}(\triangle x) / \| \triangle x \| \rar 0$ as $ \|
\triangle x \| \rar 0$.

We have $ \vl o(\triangle x+y) \vl \leq \gamma (\triangle x +y) \|
\triangle x +y \| $, where $\gamma (\triangle x +y) \rar 0$ as
$\triangle x +y \rar 0$. This is true since $diam D(x_0+\triangle
x) \rar 0$ as $\triangle x \rar 0$. The following inequalities
$$
\mid \frac{1}{\mu(D(x_0+\triangle x))} \int_{D(x_0+\triangle x)}
o(\triangle x + y) d y \mid  \leq
$$
$$
\leq \frac{1}{\mu(D(x_0+\triangle x))} \int_{D(x_0+\triangle x)}
\gamma (\triangle x +y) \| \triangle x + y \| d y \leq
$$
$$
\leq \frac{1}{\mu(D(x_0+\triangle x))} \int_{D(x_0+\triangle x)}
\gamma (\triangle x +y) (\| \triangle x \|   + \| y \|) d y .
$$
hold. Since according to the condition 2 $ \| y \| \leq k \|
\triangle x \| $ for some $k$, we get
$$
\mid \frac{1}{\mu(D(x_0+\triangle x))} \int_{D(x_0+\triangle x)}
o(\triangle x + y) d y \mid \leq \frac{\gamma ((k+1)\triangle x )
}{\mu(D(x_0+\triangle x))} \int_{D(x_0+\triangle x)}(k+1) \|
\triangle x \| d y =
$$
$$
=\gamma ((k+1)\triangle x )(k+1) \| \triangle x
\|=\tilde{o}(\triangle x).
$$
So we proved correctness of the expansion
$$
\varphi(x_0+\triangle x)= f(x_0)+(f'(x_0),\triangle
x)+\frac{1}{\mu(D(x_0+\triangle x))} \int_{D(x_0+\triangle x)}
(f'(x_0),y) d y + \tilde{o}(\triangle x).
$$
We get for $x=x_0+\triangle x$ from the regions  of constancy of
SVM $D(\cdot)$ \be
\varphi'(x)=f'(x_0)+\tilde{o}'(x-x_0).\label{firstsecsub3} \ee If
we prove that $\tilde{o}'(x-x_0) \rar 0$ as $x \rar x_0$, then it
follows from here that $\varphi'(x) \rar f'(x_0)$ as $x \rar x_0$.

We will argue in the following way. It follows from  the
definition of the infinitesimal function that $\tilde{o}'(0) =0$.
If we prove that $\tilde{o}'(\cdot)$ is a continuously
differentiable function then we will have from the above
$\tilde{o}'(\triangle x) \rar 0$ as $\triangle x \rar 0$.

Let us prove that
$$
 \varphi(x)= \frac{1}{\mu(D(x))} \int_{D(x)} f(x+y)d y,
$$
is the continuously differentiable function with respect to $x, x
\neq x_0$ if $D(\cdot)$ satisfies all requirements which can be
met easily.

Consider the function
$$
\tilde{\varphi}(x)=  \int_{D(x)} f(x+y)d y.
$$
Fix an arbitrary point $x$, an increment $\triangle x$ and
consider the difference
$$
\tilde{\varphi}(x+\triangle x) - \tilde{\varphi}(x) = \int_{D(x +
\triangle x)} f(x + \triangle x  +y)d y.- \int_{D(x )} f(x + y)d y
=
$$
$$
=\int_{D(x + \triangle x)} f(x + \triangle x  +y)d y - \int_{D(x+
\triangle x )} f(x + y)d y + \int_{D(x + \triangle x )} f(x + y)d
y -\int_{D(x )} f(x + y)d y=
$$
$$
=I_1(\triangle x)+I_2(\triangle x),
$$
where
$$
I_1(\triangle x)=\int_{D(x + \triangle x)} f(x + \triangle x  +y)d
y - \int_{D(x+ \triangle x )} f(x + y)d y,
$$
$$
I_2(\triangle x)=\int_{D(x + \triangle x )} f(x + y)d y -\int_{D(x
)} f(x + y)d y.
$$

The value $I_1(\triangle x)$ is an increment of the function
$\tilde{\varphi}(\cdot)$ in the regions of constancy of SVM
$D(\cdot)$. Only the integrand $f(\cdot)$ changes.

As it is easy to see that $I_1(\cdot)$ also depends on the
considered point $x$. The value $I_2(\triangle x)$ is a change of
the function $\tilde{\varphi}(\cdot)$ when the integrand
$f(\cdot)$ does not change and at the same time the set $D(x +
\triangle x )$, which the integration is done over, changes.

It follows from the differential qualities of the Steklov integral
(\cite{proudintegapp1},\cite{proudintegapp}) that
$$
I_1(\triangle x)=(\int_{D(x + \triangle x)} f'(x + y) d y,
\triangle x) + \overline{o}(\triangle x),
$$
where $\overline{o}(\triangle x) / \| \triangle x \| \rar 0$ as
$\triangle x \rar 0$.

We have from here
$$
I_1(\triangle x)=(\int_{D(x)} f'(x + y) d y, \triangle x) +
\hat{o}(\triangle x),
$$
where $\hat{o}(\triangle x) / \| \triangle x \| \rar 0$ as
$\triangle x \rar 0$, i.e. $I_1(\cdot)$ is continuously
differentiable at zero and \be I^{'}_1(0)=\int_{D(x)} f'(x + y) d
y. \label{firstsecsub4} \ee

Let us prove that $I_2(\cdot)$ is continuously differentiable at
zero. $I_2(\cdot)$ is also depends on $x$ like $I_1(\cdot)$.

We will prove by induction on the dimension of $\mathbb{R}^n$. It
is easy to check that for $n=1$ the function $I_2(\cdot)$ is
continuously differentiable at zero. Really, in one dimensional
case the integral
$$
\int^{b(z)}_{a(z)} f(x+y) dy
$$
is the  differentiable function with respect to $z$ with a
continuous derivative with respect to $x$ if $b(\cdot)$, $a(\cdot)
$ are continuously differentiable functions at $z=x$.

Let the statement be proved for $n=k$. We will prove it for
$n=k+1$.

Represent $I_2(\triangle x)$ in the form
$$
\int^{b(x+\Delta x)}_{a(x+\Delta x)} \theta(x+\triangle x + y_1
e_1) dy_1,
$$
where
$$
\theta(x+\Delta x+ y_1 e_1) = \int_{\hat{D}(x+\triangle x)} f(x+y)
d y_{(k)},
$$
$y=(y_1,y_2, \dots , y_k, y_{k+1})$, $e_1=(1,0,0,\dots, 0), x,
\Delta x \in \mathbb{R}^{k+1}$, $\hat{D}(x+\Delta x)$ is the
projection of $D(x+\Delta x)$ on the space $\mathbb{R}^k$
consisting from the vectors $y_{(k)}=(y_2, y_3, \dots, y_{k+1})$.

According to the induction $\theta(\cdot)$ is a continuously
differentiable function with respect to $\Delta x$. Then \be
I'_2(0)=\int^{b(x)}_{a(x)} \theta'(x+y_1 e_1) d y_1 +
\theta(x+b(x)e_1)b'(x)-\theta(x+a(x)e_1)a'(x),
\label{firstsecsub5} \ee i.e. if $a(\cdot), b(\cdot)$ are
continuously differentiable functions then $I_{2}^{'}(0)$ is a
continuous function with respect to $x$. It follows from here that
the function $\tilde \varphi(\cdot) $ is continuously
differentiable with respect to $x$ under the imposed conditions on
$a(\cdot), b(\cdot)$.

Calculate $\varphi'(\cdot)$: \be \varphi'(x)=\frac{1}{\mu(D(x))}
\tilde{\varphi}'(x)-\frac{\mu'(D(x))}{\mu^2(D(x))}
\tilde{\varphi}(x). \label{firstsecsub6} \ee

As soon as all functions in this expression are continuous
functions with respect to $x$, then $\varphi'(\cdot)$ is a
continuous function at $x \in S, x \neq x_0.$ It follows from
continuity of $\varphi'(\cdot)$ that $\tilde{o}'(\cdot)$ is a
continuous function. Consequently, $\tilde{o}'(\Delta x) \rar 0$
as $\Delta x \rar 0$.  We can make conclusion from
(\ref{firstsecsub3}) that for $x_0+\Delta x$ from the regions of
constancy of SVM $D(\cdot)$
 \be
\lim_{\Delta x \rar 0} \varphi'(x_0+\Delta x) = f'(x_0)
\label{firstsecsub7} \ee if the function $f(\cdot)$ is
differentiable at $x_0$. Notice  that the equality
(\ref{firstsecsub7}) is true for any SVM $D(\cdot)$, satisfying
the written above conditions. We get from here that for the case
when the function $f(\cdot)$ is differentiable at $x_0$ the
equality $ \p \varphi_D(x_0)=\{f'(x_0) \}$ is true.

Define SVM $ \Phi f(\cdot): \mathbb{R}^n \rar 2^{\mathbb{R}^n} $
with the images
$$
\Phi f(x_0)=\bco \, \bigcup_{D(\cdot)} \, \p \varphi_D (x_0),
$$
where the union is taken for all SVM $D(\cdot) \in \Xi$. The set
$\Phi f(x_0)$ is called {\em  the subdifferential of the first
order} of the function $f(\cdot)$ at $x_0$.

We get the following theorem from the said above
\begin{thm}
If the function $f(\cdot)$ is differentiable at $x_0$ then $\Phi
f(x_0)= \{f'(x_0) \}$. \label{thmfirstsecsub1}
\end{thm}
Let us prove some qualities of SVM $\Phi f(\cdot),$  namely, that
its images are convex compact sets.
\begin{lem}
The set $\Phi f(x_0)$ is convex and compact.
\label{lemfirstsecsub1}
\end{lem}
{\bf Proof.} The convexity is clear. Let us prove the boundedness.
We have for any $x$ from the regions of constancy of SVM
(\cite{proudintegapp1}, \cite{proudintegapp}) \be
 \varphi'(x)= \frac{1}{\mu(D(x))} \int_{D(x)} f'(x+y)d y,
\label{firstsecsub8} \ee As soon as $\| f'(x+y) \| \leq L $, then
$$
\| \varphi'(x) \| \leq L.
$$
Consequently, $\Phi f(x_0)$ is a convex bounded set. Let us prove
the closure.

Let  $v_i \in \p \varphi_{D_i}(x_0)$ and $v_i \rar v$ for SVM
$D_i(\cdot) \in \Xi$. Let us prove that $v \in \Phi (x_0)$.
Construct from SVM $D_i(\cdot)$ a new SVM $D(\cdot)$ from the set
$ \Xi $. For this it is sufficiently to consider a subsequence $\{
i_k \} \subset \{ i \} $ such that $j \in \{ i_k \}$ and for a
point $x_j,$ corresponding $v_j,$ the equality $D(x_j)=D_j(x_j)$
and the inequality
$$
\|  \frac{1}{\mu(D_j(x_j))} \int_{D_j(x_j)} f'(x_j+y)d y - v \|
\leq \eps_j,
$$
where $\eps_j \xrightarrow[j]{} +0,$ were correct in some
neighborhood of $x_j$ from the regions of constancy of SVM
$D_j(\cdot)$.

It is true for the constructed SVM $D(\cdot)$ at the points $x_j$
from the regions of its constancy
$$
\lim_{x_j \rar x_0} \varphi'(x_j) = \lim_{x_j \rar x_0}
\frac{1}{\mu(D(x_j))} \int_{D(x_j)} f'(x_j+y)d y   = v,
$$
i.e. $v \in \Phi f(x_0)$. Consequently, the set $\Phi f(x_0)$ is
closed. The compactness follows from the boundedness and closure.
The lemma is proved. $\Box$

A set $\eta(x_0)$ of smooth curves in $\mathbb{R}^n$ was defined
in \cite{lowapp2}, \cite{lowapp2a} to analyse the differential
qualities of the function $f(\cdot)$.
\begin{defi} $ \eta (x_0)$  is the set of the smooth curves $r(x_0,\al,
g) = x_0 + \al g+ o_r(\al)$  where $g \in S^{n-1}_1(0)=\{ v \in
\mathbb{R}^n \vl \| v \| =1 \}$ and the function
$o(\cdot):[0,\al_0] \rar R^n, \al_0>0$ satisfies the next
conditions

\nid 1) $o_r(\al)/(\al) \rar +0$ uniformly in $r(\cdot)$ as $\al
\rar +0$

\nid 2) there is the continuous derivative $o_r'(\cdot)$ and its
norm is bounded for all $r$ in the following sense:   $c \; <
\infty $ exists such that
$$
     \sup_{\tau \in (0,\al_0)} \parallel o_r'(\tau) \parallel \leq c
$$
3) the derivative $\grad f(r(\cdot))$ exists almost everywhere
(a.e.).
 along the curve  $r(x_0,\cdot ,g)$ \\
\label{defi1321}
\end{defi}
\begin{rem}
According to the property 3 of the definition the set $\eta(x_0)$
depends on choosing $f(.).$ \label{rem321}
\end{rem}
Consider for some $g \in S^{n-1}_1(0)$ a curve $r(\cdot) \in \eta
(x_0),$ that is defined on the segment $[0,\al_0]$. Take any
sequence $\{\al_k\}, \al_k \rar +0,$ as $k \rar \infty $ and
consider the average integral limit values of the gradients $\grad
f(r(\cdot))$ along such curves $r(\cdot)$
$$
 \avin{0}{{\al_k}}{{\al_k}}{x_0}{g}{f}.
$$
The limit value of these averages as $k \rar \infty$ contains an
important information about behavior of the function $f(\cdot)$
near the point $x_0$ in the direction $g$.

Introduce the sets
$$
E f(x_0) =  \{ v \in R^n : \exs {\al_{k}}, \al_{k} \rar +0, ( \exs
\,g \in S^{n-1}_1 (0) ),
$$
$$(\exs r(x_0,\cdot ,g)\, \in \, \eta(x_0)) ,
 v = \lim_{\al_k \rar +0} \al_k^{-1} \; \int^{\al_k}_0 \,
 \grad f(r(x_0,\tau,g))d\tau \; \}
$$
and
$$
Df(x_0)= \mbox{co} \,\, Ef(x_0),
$$
where  the integral is the Lebesgue integral
\cite{kolmogorovfomin}.

Let us prove the following theorem.
\begin{thm}
It is true
$$
\Phi f(x_0)=Df(x_0).
$$
\label{thmfirstsecsub2}
\end{thm}
{\bf Proof.} The derivative of the function $\varphi(\cdot)$ at
points $x$ from the regions of constancy of SVM $D(\cdot) \in \Xi$
is calculated according to the formula (\ref{firstsecsub8}).
Rewrite (\ref{firstsecsub8}) in the form of the integral sums:
$$
\varphi'(x)=\frac{1}{\mu(D(x))} \lim_{N \rar \infty}
\sum_{i=1}^{N} f'(x+y_i) \mu(\Delta D_i)= \lim_{N \rar \infty}
\sum_{i=1}^{N} f'(x+y_i) \frac{\mu(\Delta D_i)}{ \mu(D(x))}=
$$
\be =\lim_{N \rar \infty} \sum_{i=1}^{N} f'(x+y_i) \beta_i,
\label{firstsecsub9} \ee where $y_i \in \Delta D_i$, $\beta_i
=\frac{\mu(\Delta D_i)}{ \mu(D(x))}$, $0 \leq \beta_i \leq 1$,
$\sum_{i=1}^{N} \beta_i=1$, $D(x)=\bigcup_{i=1}^{N} \Delta D_i$,
$\mu(D(x))=\sum_{i=1}^{N} \mu(\Delta D_i)$. we can make conclusion
from here that $\varphi'(x)$ is the convex envelope of the vectors
$ f'(x+y_i).$

Divide the region of integration $x+D(x)$ into sectors (cones),
not having inner common points, with common vertex $x_0$. Take in
each $i-$ th cone (sector) a curve $r(x_0,\cdot, g_i) \in
\eta(x_0).$ The integral
$$
\al^{-1} \; \int^{\al}_0 \, \grad f(r(x_0,\tau,g_i))d\tau
$$
can be represented as the limit of the convex envelope of the
gradients of the function $f(\cdot)$, calculated along the curve
$r(x_0,\cdot, g_i)$. Really, \be \al^{-1} \; \int^{\al}_0 \, \grad
f(r(x_0,\tau,g_i))d\tau= \lim_{N \rar \infty} \sum_{j=1}^{N}
f'(r(x_0,\tau_j,g_i) ) \gamma_j, \label{firstsecsub10}\ee where
$\gamma_j=\frac{\Delta \tau_j}{\al}$, $\Delta \tau_j $ is a
segment of integration with respect to $\tau$ along the curve
$r(x_0,\tau_j,g_i)$, $\bigcup_{i=1}^{N} \Delta \tau_i = [0,\al], $
$\sum_{j=1}^{N}\gamma_j=1 $, $0\leq \gamma_j \leq 1$.

Place each curve $r(x_0,\tau,g_i)$ into a curved cylindrical body
$\Delta S_i$ with any small measure $\mu(\Delta S_i)$. Divide
$\Delta S_i$ with the help of planes $\pi_{ij}$, normal to the
central axis of $\Delta S_i$, into pieces $\Delta S_{ij}$ with the
measures $\mu(\Delta S_{ij})$. Denote the measure of intersection
and $\pi_{ij} \cap \Delta S_i $ by $\Delta \varphi_i $. Then $\al
\cdot \Delta \varphi_i \simeq   \mu(\Delta S_i),$ $\Delta \tau_j
\cdot \Delta \varphi_i \simeq \mu(\Delta S_{ij}).$ Consequently,
$\gamma_j= \frac{\Delta \tau_j \cdot \Delta \varphi_i}{\al \cdot
\Delta \varphi_i} \simeq \frac{\mu(\Delta S_{ij})}{\mu(\Delta
S_i)}$.

It is obvious, that the convex envelope on $g_i$ of
(\ref{firstsecsub10}) is the private case of the convex envelope
(\ref{firstsecsub9}) for the reason of arbitrary choice of SVM
$D(\cdot) \in \Xi$ and any small measure of intersection of the
cylindrical bodies $\Delta S_i$. It follows from here \be Df(x_0)
\subset \Phi f(x_0). \label{firstsecsub11} \ee

But from the other side the integral sum (\ref{firstsecsub9}) can
be considered as the private case of the integral sum
(\ref{firstsecsub10}), if to divide the set $x+D(x)$ into sectors
$\Delta S_i$, not having inner points,  with the common vertex
$x_0$, as soon as we are free in choosing of way of dividing into
$\Delta S_{ij}.$

We will choose in each sector (cone) a curve $r(x_0,\cdot, g_i)
\in \eta(x_0).$ Place $i-$th sector (cone) into a cylindrical body
$\Delta S_i$ in such way that the measure of $i$-the sector (cone)
were equal to half of the measure of $\Delta S_i$. But integration
over the cylindrical body $\Delta S_i$ is equal to integration
over $i-$th the sector (cone) two times.

Divide each cylindrical body $\Delta S_i$ by the planes
$\pi_{ij}$, normal to the central axis of $\Delta S_i$, into the
pieces $\Delta S_{ij}$. Take in $\Delta S_{ij}$ the points
$x+y_{ij} = r(x_0,\tau_j,g_i) \in \Delta S_{ij}$.

As s result, the convex envelope of the gradients of the function
$f(\cdot)$ at $x+y_{ij}$ has the form \be \sum_{j=1}^{N}
f'(x+y_{ij}) \beta_{ij}, \label{firstsecsub12} \ee  where
$\beta_{ij} =\frac{\mu(\Delta S_{ij})}{ \mu(\Delta S_{i}))}$, $0
\leq \beta_{ij} \leq 1$, $\sum_{j=1}^{N} \beta_{ij}=1$, $\Delta
S_i=\bigcup_{j=1}^{N} \Delta S_{ij}$, $\mu(\Delta
S_i)=\sum_{j=1}^{N} \mu(\Delta S_{ij})$. Instead of taking of the
convex envelope of the gradients of $f(\cdot)$ at the points from
the whole set $x+D(x)$, we, at first, take the convex envelope of
the gradients of $f(\cdot)$ in each cylindrical body $ \Delta S_i
$. As a result, we get the sum (\ref{firstsecsub12}) which is
approximately equal to the sum (\ref{firstsecsub10}). The bigger
$N$ and  smaller sectors (cones), the more precise equality for
the sums will be. To get the sum (\ref{firstsecsub9}) it is
necessarily to take the convex envelope on $i$ of the sums
(\ref{firstsecsub12}), as soon as the convex envelope of the
convex envelope is the convex envelope again. Finally, we get that
the sum (\ref{firstsecsub9}) can be obtained as the convex
envelope on $i$ of the sums (\ref{firstsecsub10}). The said is
correct for any SVM $D(\cdot) \in \Xi$. From here we have \be \Phi
f(x_0) \subset D f(x_0). \label{firstsecsub13} \ee The statement
of the theorem follows from (\ref{firstsecsub11}) and
(\ref{firstsecsub13}). $\Box$

The proved theorem agrees with the earlier proved Theorem
\ref{thmfirstsecsub1} for the case when $f(\cdot)$ is
differentiable at $x_0$  and also with Lemma \ref{lemfirstsecsub1}
for any Lipschitz function $f(\cdot)$ as soon as for these cases
 $\Phi f(x_0)=D f(x_0)=\{ f'(x_0) \}$ and $D f(x_0)$
is a convex compact set (see \cite{lowapp2}).

\section{\bf  The subdifferential of the second order}

We will consider the function $\psi(\cdot):\mathbb{R}^n \rar
\mathbb{R}$
$$
 \psi(x)= \frac{1}{\mu(D(x))} \int_{D(x)} \varphi(x+y)d y,
$$
to construct the subdifferential of the second order where
$\varphi(\cdot)$ was defined before for SVM $D(\cdot) \in \Xi$.
The function $\psi(\cdot)$ depends on the chosen SVM $D(\cdot) \in
\Xi$.

If SVM $D(\cdot)$ is constant then, how it was proved in
\cite{proudintegapp}, $\psi(\cdot)$ is a twice differentiable
function. Our goal is to define a set consisting of the
generalized matrices at the point $x_0$ for the Lipschitz function
$f(\cdot)$ in such way that for any twice differentiable function
the set of the generalized matrices would coincide with the matrix
of the second mixed derivatives of this function. We have the
similar situation for a differentiable function and the
subdifferential of the first order because the last one coincides
with the derivative of this function.

Define for SVM $D(\cdot) \in \Xi$ similar to that, how it was done
above, the set
$$
\p \psi_D(x_0)=\bco \{ v \in \mathbb{R}^n \vl v= \lim_{x_i \rar
x_0} \psi'(x_i) \},
$$
where the points $x_i$ are taken from the regions of constancy of
SVM $D(\cdot)$. Similar to the proof of Lemma
\ref{lemfirstsecsub0} we can prove that $\p \psi_D(x_0)$ is a
convex compact set.

Introduce SVM $ \Psi f(\cdot): \mathbb{R}^n \rar 2^{\mathbb{R}^n}
$ with the images
$$
\Psi f(x_0)=\bco \, \bigcup_{D(\cdot)} \, \p \psi_D (x_0),
$$
where the union is taken over all SVM $D(\cdot) \in \Xi$. How it
follows from the theorem, proved below, the set $\Psi f(x_0)$ can
be called the subdifferential of the first order of the function
$f(\cdot)$  at $x_0$ as well.

Similar to the proof of Lemma \ref{lemfirstsecsub1} we can prove
that $\Psi f(x_0)$ is a convex compact set. It appears that it
coincides with $Df(x_0)$.
\begin{thm}
The equality
$$
\Psi f(x_0)=Df(x_0).
$$
is true.
\end{thm}
{\bf Proof.} We will prove in two steps similar to the proof of
Theorem \ref{thmfirstsecsub2}. The set $\p \psi_D (x_0)$ consists
of the limit values of vectors, equaled to the convex envelope of
the gradients of the function $\varphi(\cdot)$ at the points $z$
of $x+D(x)$, where the points $x$ are taken from  the regions of
constancy of SVM $D(\cdot)$, and $x \rar x_0$. But the gradient of
the function $\varphi(\cdot)$ at any point $z \in x+D(x)$ is equal
to the convex envelope of the gradients of the function $f(\cdot)$
at the points of $x+2D(x)$, where these gradients exist, and the
points $x$ are taken from the regions of constancy of SVM
$D(\cdot)$.

It is known that the convex envelope of  vectors $\{a\}$, when
each vector $a$ is the convex envelope of vectors $\{b\}$, is the
convex envelope of vectors $\{b\}$. Therefor, the vectors of the
set $\p \psi_D (x_0)$ are the convex envelope of the gradients of
the function $f(\cdot)$ at the points of $x+2D(x)$, where these
gradients exist, and $x$ are taken from the regions of constancy
of SVM $D(\cdot)$ as $x \rar x_0$.

We have already proved Theorem \ref{thmfirstsecsub2} based on the
fact that the gradients of the function $\varphi(\cdot)$ are equal
to the convex envelope of the gradients of the function $f(\cdot)$
at the points $z \in x+D(x)$ where these gradients exist. The
equality $\Phi f(x_0)=Df(x_0)$ has been already proved.

As soon as the gradient of the function $\psi(\cdot)$ at $x$ is
the convex envelope of the gradients of the function $f(\cdot)$ at
the points of the set $x+2D(x)$, where they exist, then repeating
the arguments of Theorem \ref{thmfirstsecsub2}, we will get the
equality $ \Psi f(x_0)=Df(x_0).$ The theorem is proved. $\Box$

The next step is to give a definition of the generalized matrices
of the function $f(\cdot)$ at $x_0$ and also a definition of the
subdifferential of the second order $\Psi^2 f(x_0)$, consisting of
the generalized matrices.

Introduce a set of the matrices
$$
\p^2 \psi_D(x_0)=\bco \{ A \in \mathbb{R}^{n \times n} \vl A=
\lim_{x_i \rar x_0} \psi''(x_i) \},
$$
where the points $x_i$ belong to the regions of constancy of SVM
$D(\cdot) \in \Xi $.
\begin{lem}
$\p^2 \psi_D(x_0)$ is a closed convex  set.
\label{lemfirstsecsub2}
\end{lem}
{\bf Proof}. The convexity is clear. Let us prove the closure. Let
$A_i \in \partial^2 \psi_{D}(x_0)$ and $A_i \rar A$ for SVM
$D(\cdot) \in \Xi$. Let us prove, that $A \in \p^2 \psi_{D}(x_0)$.
Consider a subsequence $\{ i_k \} \subset \{ i \} $ such that for
$j \in \{ i_k \}$ and the points $\{ x_j \}$ from the regions of
constancy of SVM $D(\cdot)$, for that the limit of  $\psi''(x_j)$
is equal to $A_j$, the inequality
$$
\|  \frac{1}{\mu(D(x_j))} \int_{D(x_j)} \varphi''(x_j+y)d y - A \|
\leq \eps_j,
$$
would be true, where $\eps_j \xrightarrow[j]{} +0. $ It is true
$$
\lim_{x_j \rar x_0} \psi''(x_j) = \lim_{x_j \rar x_0}
\frac{1}{\mu(D(x_j))} \int_{D(x_j)} \varphi''(x_j+y)d y   = A,
$$
for  the points $x_j$ from the regions of constancy of SVM
$D(\cdot)$ i.e. $A \in \p^2 \psi_{D}(x_0)$. Consequently, the set
$\p^2 \psi_{D}(x_0)$ is closed. The lemma is proved. $\Box$

Define SVM $ \Psi^2 f(\cdot): \mathbb{R}^n \rar 2^{\mathbb{R}^{n
\times n}} $ with the images
$$
\Psi^2 f(x_0)=\bco \, \bigcup_{D(\cdot)} \, \p^2 \psi_D (x_0),
$$
where the union is taken over all SVM $D(\cdot) \in \Xi$. The set
$\Psi^2 f(x_0)$ is {\em called the subdifferential of the second
order} of the function $f(\cdot)$ at the point $x_0$.

Let us prove some qualities of this set.
\begin{lem}
$\Psi^2 f(x_0) $ is a convex closed set.
\end{lem}
{\bf Proof}. The convexity is clear. Let us prove the closure. We
will prove it by the same method how we did it in Lemma
\ref{lemfirstsecsub1}. Let  $A_i \in \p^2 \psi_{D_i}(x_0)$ and
$A_i \rar A$ for SVM $D_i(\cdot) \in \Xi$. We will prove that $A
\in \Psi^2 f(x_0)$. Compose from SVM $D_i(\cdot)$ a new SVM
$D(\cdot) \in \Xi$.  It is sufficiently for this to consider a
subsequence $\{ i_k \} \subset \{ i \} $ such that for $j \in \{
i_k \}$ and for the points $\{ x_j \} ,$ corresponding to the
matrix $A_j,$ in some their surroundings from the regions of
constancy of SVM $D_j(\cdot)$ the equalities $D(x_j)=D_j(x_j)$ and
also the inequality
$$
\|  \frac{1}{\mu(D_j(x_j))} \int_{D_j(x_j)} \varphi''(x_j+y)d y -
A \| \leq \eps_j,
$$
would be correct, where $\eps_j \xrightarrow[j]{} +0. $

The equality
$$
\lim_{x_j \rar x_0} \psi''(x_j) = \lim_{x_j \rar x_0}
\frac{1}{\mu(D(x_j))} \int_{D(x_j)} \varphi''(x_j+y)d y   = A,
$$
is correct for  SVM $D(\cdot)$ at $x_j$ from the regions of
constancy i.e. $A \in \Psi^2 f(x_0)$. Consequently, the set
$\Psi^2 f(x_0)$ is closed. The lemma is proved.$\Box$
\begin{rem}
Remark, that the sets  $\p^2 \psi_D(x_0)$ and $\Psi^2 f(x_0)$ may
be unbounded without some additional assumptions for the function
$f(\cdot)$.
\end{rem}
Consider the case when $f(\cdot)$ is twice differentiable at
$x_0$. The equality
$$
f(x_0+\Delta x)=f(x_0)+(f'(x_0), \Delta x)+\frac{1}{2}(f''(x_0)
\Delta x, \Delta x)+o( \| \Delta x \|^2),.
$$
is true where $o( \| \Delta x \|^2) /  \| \Delta x \|^2  \rar 0$
as $\Delta x \rar 0$. Find the answer for the question: what is
the set $\Psi^2 f(x_0)$ in this case?

Write down an expression for the function $\varphi(\cdot)$. As
soon as the matrix  $f''(x_0)$ is symmetric, we will have
$$
\varphi(x_0+\triangle x)= f(x_0)+(f'(x_0),\triangle x)+
\frac{1}{2}(f''(x_0) \Delta x, \Delta x)+
$$
$$
+\frac{1}{\mu(D(x_0+\triangle x))} \int\limits_{D(x_0+\triangle
x)} (f'(x_0),y) d y + \frac{1}{\mu(D(x_0+\triangle x))}
\int\limits_{D(x_0+\triangle x)} (f''(x_0) \Delta x,y) d y +
$$
$$
+\frac{1}{2\mu(D(x_0+\triangle x))} \int\limits_{D(x_0+\triangle
x)} (f''(x_0) y,y) d y+ \frac{1}{\mu(D(x_0+\triangle x))}
\int\limits_{D(x_0+\triangle x)}{o}( \| \triangle x +y \|^2)dy.
$$
Denote
\begin{eqnarray}
& & \Theta (\Delta x)=\frac{1}{\mu(D(x_0+\triangle x))}
\int\limits_{D(x_0+\triangle x)} (f'(x_0),y) d y+
\frac{1}{\mu(D(x_0+\triangle x))} \int\limits_{D(x_0+\triangle x)}
(f''(x_0) \Delta x,y) d y +{}\nonumber\\
&&+\frac{1}{2\mu(D(x_0+\triangle x))} \int\limits_{D(x_0+\triangle
x)} (f''(x_0) y,y) d y.{}\nonumber
\end{eqnarray}
It is obvious that $ \Theta(\cdot)$ is a linear function with
respect to $\Delta x$, when $x+\Delta x$ belongs to the regions of
constancy of  SVM $D(\cdot)$. In this case the function
$\varphi(\cdot)$ has the form
$$
\varphi(x_0+\triangle x)= f(x_0)+(f'(x_0),\triangle x)+
\frac{1}{2}(f''(x_0) \Delta x, \Delta x)+
$$
$$
+\Theta(\Delta x)+ \frac{1}{\mu(D(x_0+\triangle x))}
\int\limits_{D(x_0+\triangle x)}{o}( \| \triangle x +y \|^2)dy.
$$
The function $\psi(\cdot)$ can be written in the form
$$
\psi(x_0+\triangle x)= f(x_0)+(f'(x_0),\triangle x)+
\frac{1}{2}(f''(x_0) \Delta x, \Delta x)+\tilde\Theta(\Delta x)+
$$
\be + \frac{1}{\mu(D(x_0+\triangle x))}
\int\limits_{D(x_0+\triangle x)} \left(
\frac{1}{\mu(D(x_0+\triangle x+z))} \int\limits_{D(x_0+\triangle
x+z)}{o}( \| \triangle x +y+z \|^2) dy \right)   dz,
\label{firstsecsub19} \ee where $\tilde\Theta(\cdot)$ is the
linear function with respect to  $\Delta x $ when  $x+\Delta x$
from the regions of constancy of SVM $D(\cdot)$.

Let us prove, that for SVM $D(\cdot) \in \Xi $
\begin{eqnarray}
& &\frac{1}{\mu(D(x_0+\triangle x))} \int\limits_{D(x_0+\triangle
x)} \left( \frac{1}{\mu(D(x_0+\triangle x+z))}
\int\limits_{D(x_0+\triangle x+z)}{o}( \| \triangle x +y+z \|^2)
dy \right)   dz=  {}\nonumber \\
& &  =\tilde{o}(\| \triangle x \|^2),  {}\nonumber
\end{eqnarray}
where $\tilde{o}(\| \triangle x \|^2) / \| \triangle x \|^2 \rar
0$ as $ \| \triangle x \| \rar 0$.

According to the qualities of SVM $D(\cdot) \in \Xi:$   $\| y \|
\leq k \| \Delta x +z \| $ and $o(\|\Delta x+z+y \|^2) \leq \gamma
\|\Delta x+z+y \|^2 $, where $ \gamma = \gamma(\|\Delta x+z+y \|)
\rar 0 $, as $ \|\Delta x+z+y \| \rar 0$, then
$$
o(\|\Delta x+z+y \|^2) \leq  2 \gamma (\|\Delta x+z \|^2 + \| y
\|^2) \leq 2\gamma ( \|\Delta x+z \|^2 + k^2 \|\Delta x+z \|^2) =
$$
$$
= 2 \gamma (1+k^2)\|\Delta x+z \|^2.
$$
We have used the inequality $2ab \leq a^2+b^2 $. From here
$$
(a+b)^2 = a^2+b^2+2ab \leq 2(a^2 + b^2).
$$
According to the qualities of SVM $D(\cdot) \in \Xi$ we have $\| z
\| \leq k \| \Delta x \|$. Consequently,
$$
o(\|\Delta x+z+y \|^2) \leq 4 \gamma (1+k^2)(\|\Delta x \|^2  + \|
z \|^2) \leq 4 \gamma (1+k^2)^2 \|\Delta x \|^2.
$$
It follows from here that
\begin{eqnarray}
&&\frac{1}{\mu(D(x_0+\triangle x))} \int\limits_{D(x_0+\triangle
x)} \left( \frac{1}{\mu(D(x_0+\triangle x+z))}
\int\limits_{D(x_0+\triangle x+z)}{o}( \| \triangle x +y+z \|^2)
dy \right)   dz \leq {}\nonumber \\
&&\leq \gamma  (1+k^2)^2 \|\Delta x \|^2 = \tilde{o}(\| \triangle
x \|^2), {}\nonumber
\end{eqnarray}
i.e. $ \tilde{o}(\| \triangle x \|^2) /  \| \triangle x \|^2 \rar
0$, as $\Delta x \rar 0$, and $\gamma \rar 0$.

It follows from the formula (\ref{firstsecsub19}) that the
equality \be \psi''(x)= f''(x_0)+ \tilde{o}''(\|  x - x_0 \|^2)
\label{firstsecsub20} \ee is correct for the twice differentiable
function $f(\cdot)$ at $x_0$ and $x \neq x_0$ from the regions of
constancy of SVM $D(\cdot)$.

We have  from the definition of the infinitesimal function of the
second order that  $ \tilde{o}(\cdot) $ has the first and second
derivatives, equaled to zero, at the point $x_0$. Show that $
\tilde{o}(\cdot)$ is a twice differentiable function. It follows
from here $\tilde{o}'(\triangle x)  \rar 0$  and $\tilde{o}''(\|
\triangle x \|^2) \rar 0$ as $\triangle x \rar 0$.

It is not difficult to prove, that $\psi(\cdot)$ is a twice
differentiable function for $ x \neq x_0$ if the boundary of the
set $D(x)$ is given by twice differentiable functions with respect
to $x$.

As soon as the function $\tilde \Theta(\cdot)$ is expressed in
terms of integrals of twice continuously  differentiable functions
with respect to $\triangle x$, then $\tilde \Theta(\cdot)$ is a
twice continuously  differentiable function.

Consequently, $\tilde{o}'(\triangle x)  \rar 0$ and
$\tilde{o}''(\| \triangle x \|^2) \rar 0$ as $\triangle x \rar 0$.
It follows from (\ref{firstsecsub20}) that $\psi''(x_0+\triangle
x) \rar f''(x_0)$ when $x+\triangle x$ from the regions of
constancy of SVM $D(\cdot)$ as $\triangle x \rar 0$. The following
theorem is proved.
\begin{thm}
If $f(\cdot)$ is a twice differentiable function at $x_0$, then
$$\Psi^2 f(x_0) = \{ f'' (x_0) \} $$.
\label{firstsecondsubdifftheorem32}
\end{thm}

\vspace{0.5cm}
\section{\bf Application of the subdifferentials of the first and
second orders}

The necessary condition of optimality can be written in different
ways. Let us write down one of them.

Let  $\Omega \subset \mathbb{R}^n$  be a convex compact set.
Define for any point $x_0 \in \Omega$ a cone of the tangent
directions
$$
K(x_0, \Omega) = \{ g \in \mathbb{R}^n \vl \exists \beta_0 >0,
\exists r(x_0,\al,g)=x_0+\al g +o(\al) \in \eta(x_0),
$$
\be o(\al)/\al \rar_{\al \rar +0} +0, \, r(x_0,\al,g) \in \Omega
\,\,\,\,\, \forall \al \in [0,\beta_0] \}. \label{1334} \ee Form a
set of  limit vectors
$$
A(x_0)=\bco \{v(g) \in \mathbb{R}^n \vl \exists \{ \al_k \}, \al_k
\rar_k +0, \exists g \in K(x_0,\Omega), \exists r(x_0,\cdot,g) \in
\eta(x_0) :
$$
$$
v(g)=\lim_{\al_k \rar +0} \al_k^{-1} \int^{\al_k}_0 \grad
f(r(x_0,\tau,g)) d \tau \},
$$
where $r(x_0.\al,g) \in \Omega$ for small $\al.$
\begin{lem}
For  $x_* \in \Omega$ to be a minimum point of $f(\cdot)$ on the
set $\Omega$ it is necessarily that\be \max_{v \in D f(x_*)} (v,g)
\meq 0 \;\; \forall g \in K(x_*,\Omega). \label{1335} \ee
\label{lem1331}
\end{lem}
{\bf Proof}. For any $v(g) \in A(x_*)$ there exists a sequence
$\{\al_k\}, \al_k \rar_k +0, g \in K(x_*, \Omega)$ and
$r(x_*,\cdot,g) \in \eta(x_0),$ that
$$
f(r(x_*,\al_k,g))=f(x_*)+\al_k(v(g),g)+o(\al_k) \;\; \forall k,
$$
where $o(\al_k) \rar 0$ as $\al_k \rar_k 0.$ If $x_*$ is a minimum
point of the function $f(\cdot)$ on $\Omega$, then
$$    (v(g),g) \meq 0     \;\; \forall g \in K(x_0,\Omega).   $$
From here
$$
0 \leq (v(g),g) \leq \max_{v \in A(x_*)} (v,g) \leq \max_{v \in D
f(x_*)} \, (v,g) \;\; \forall g \in K(x_*,\Omega).
$$
The lemma is proved. $\triangle$
\begin{rem}
It follows from the proof of Lemma \ref{lem1331}, that the
necessary condition of the minimum  of the function $f(\cdot)$ at
$x_*$ can be written in the form
$$
        \max_{v \in A(x_*)}\, (v,g) \meq 0 \;\; \forall g \in K(x_0,\Omega).
$$
\label{rem1331}
\end{rem}
\begin{rem}
If the function $f(\cdot)$ is convex, then (\ref{1335}) is the
sufficient condition for the minimum of $f(\cdot)$ at $x_*$, since
in this case according to the results, published in
\cite{lowapp2},
$$
Df(x_*)= \p f(x_*) .
$$
\end{rem}
\begin{cor}
For a point $x_*$ to be the minimum of the function $f(\cdot)$, it
is necessarily the condition
$$
     \min_{v \in Df(x_*)} \; (v,g) \leq 0 \;\; \forall g \in K(x_*,\Omega),
$$
or
$$
    \min_{v \in A(x_*)} (v,g) \leq 0 \;\; \forall g \in K(x_*,\Omega)
$$
was satisfied.
\end{cor}
Denote through $K^+(x_0,\Omega)$ the conjugate cone of the cone
$K(x_0,\Omega),$ that by definition is \be K^+(x_0,\Omega)=\{ w
\in \mathbb{R}^n \vl (v,w) \meq 0 \;\; \forall v \in K(x_0,\Omega)
\}. \label{1336} \ee
\begin{thm}
For a point $x_* \in \Omega $ to be the minimum of the function
$f(\cdot)$ on the set $\Omega$, it is necessarily that \be Df(x_*)
\cap K^+(x_*,\Omega) \neq \O. \label{1337} \ee \label{thm1333}
\end{thm}
{\bf Proof.} Let \be Df(x_*) \cap K^+(x_*,\Omega) = \O.
\label{1338} \ee Let us use the following theorem.
\begin{thm} (\cite{demvas})
Let $K \subset \mathbb{R}^n$ be a closed convex cone and $G
\subset \mathbb{R}^n$ be a convex compact. In order for sets $K$
and $G$ have no common points i.e.
$$         K \cap  G = \O,               $$
it is necessarily  and sufficiently that a vector $w_0 \in K^+$
existed that
$$
         \max_{x \in G} \, (w_0,x) < 0.
$$
\end{thm}
In our case (\ref{1338}) means that a vector $\bar{g} \in
K^{++}(x_*,\Omega)=K(x_*,\Omega)$  exists such that \be \max_{w
\in Df(x_*)} \, (w,\bar{g}) <0. \label{1339} \ee For the vector
$\bar{g} \in K(x_*,\Omega)$ there exist a vector $v(\bar{g}) \in
Df(x_*),$ a curve $r(x_*,\cdot,\bar{g}) \in \eta(x_*)$ and a
sequence $\{ \al_k \}, \al_k \rar_k +0,$ that $$
v(\bar{g})=\lim_{\al_k \rar +0} \al_k^{-1} \int^{\al_k}_0 \grad
f(r(x_*,\tau,\bar{g})) d \tau.
$$
The last one means that the expansion
$$
f(r(x_*,\al_k,\bar{g}))=f(x_*)+\al_k(v(\bar{g}),\bar{g})+o(\al_k),
$$
is true where $o(\al_k)/\al_k  \rar_k 0.$ From (\ref{1339}) it
follows that \be (v(\bar{g}),\bar{g}) \leq \max_{w \in Df(x_*)} \,
(w,\bar{g}) <0. \label{13310} \ee From here
$$
\frac{f(r(x_*,\al_k,\bar{g}))-f(x_*)}{\al_k} =
(v(\bar{g}),\bar{g})+ \frac{o(\al_k)}{\al_k}.
$$
Taking into account the quality of the function $o(\cdot)$ and
(\ref{13310}) we have  for small $\al_k>0$
$$
         f(r(x_*,\al_k,\bar{g})) < f(x_*).
$$
We got the contradiction with the statement that $x_*$ is the
minimum. The theorem is proved. $\triangle$
\begin{rem}
How it was mentioned in Remark \ref{rem1331}, we can consider the
set $A(x_*)$ instead of $Df(x_*)$. Then the condition (\ref{1335})
can be rewritten in the form
$$
A(x_*) \cap K^+(x_*,\Omega) \neq \O.
$$
\end{rem}
\begin{cor}
The necessary condition of optimality in $\mathbb{R}^n$ is
$$    0 \in Df(x_*).         $$
\label{corfirstsecsub2}
\end{cor}
{\bf Proof.} We will prove the statement under condition that $
x_*$ is the minimum. The statement for the maximum can be proved
in analogous way.

If $\Omega = \mathbb{R}^n,$ then $K^+(x_*,\Omega)=\{0\}.$
Consequently, the conclusion $0 \in Df(x_*)$ follows from Theorem
\ref{thm1333}.  $\triangle$

We are looking for the necessary and sufficient conditions of
optimality.

Consider the directional derivative of the function $f(\cdot) $ in
a direction $g \in S^{n-1}_1(0)$ at $x_0$ which, by definition, is
$$
\frac{\p^{\downarrow} f(x_0)}{\p g}=\underline\lim_{ \,\,\al_k
\rar +0} \frac{f(x_0+\al_k g)-f(x_0)}{\al_k}= \max_{v \in B} (v,g)
- \max_{w \in A} (w,g).
$$
The necessary condition of the minimum at $x_0$ of $f(\cdot)$  is
$A \incl B$.

If $A$ and $B$ have a common point $v$ on the boundaries, the
unity over which we denote by $\Upsilon$, then there is a set $G$
of some suspicious directions $g \in S^{n-1}_1(0)$  for the
extremum, where $S^{n-1}_1(0)=\{v \in \mathbb{R}^n \vl \| v \|=1
\}$ is the unit sphere with center at $0$. The set $G$ is the
unity over $v \in \Upsilon$ of intersections of the unit sphere
$S^{n-1}_1(0)$ and the normal cones to $A$ and $B$, constructed
for $v \in \Upsilon$.

The set $G$ can be covered by cones $K(v)$, $ v \in \Upsilon,$
with the common vertex at the point $0$. As soon as the function
$\psi(\cdot) $ for any SVM $D(\cdot) \in \Xi$ is the limit of the
convex envelopes of values of the function $f(\cdot)$ at the
points from the set $D(x)$, defined in a small neighborhood of the
point $x_0$, then the function $\psi(\cdot) $ will behave itself
in the suspicious direction $g$  like the function $f(\cdot)$. It
means that if $f(x_0 + \al_k g) > f(x_0)$ for small $\al_k
>0$, then some functions $\psi(\cdot)$ satisfy the same inequality
$\psi(x_0 + \al_k g) > \psi(x_0)$. The same is true for the case
when for small $\al_k >0$, $\al_k \rar_k  +0,$ the inequality
$f(x_0 + \al_k g) < f(x_0)$ is true.

It follows from the said above and Theorems \ref{thmfirstsecsub1},
\ref{firstsecondsubdifftheorem32}  that the following theorem can
be a sufficient condition of optimality.
\begin{thm}
If the necessary condition for the minimum of $f(\cdot)$ at $x_0$
is true  and there exists $\beta(g) >0$ for all suspicious
directions $g \in G$ that the inequality
$$
(Ag,g) \meq \beta(g) \| g  \|^2 \,\,\,\,\,\, \forall A \in \Psi^2
f(x_0),
$$
is true, then $x_0$ is the minimum of $f(\cdot)$.
\label{thmfirstsecsub5}
\end{thm}

Theorem \ref{thmfirstsecsub5} requires that all matrices $A \in
\Psi^2 f(x_0)$ were positive definite. This condition is too
heavy. Indeed, it is sufficiently to demand positive definiteness
of $A \in \p^2 \psi_D f(x_0)$ for such SVM $D(\cdot) \in \Xi$ that
cover the sets $D(x) \cap  K(v)$, $v \in \Upsilon$, as precisely
as possible for all points $x=x_0+ \alpha g$, $\al \in
\mathbb{R}^+ $, $g \in K(v)$.

Define the set of SVM $D(\cdot) \in \Xi$, that satisfy the said
above, by $\Im \incl \Xi$. Then the less heavy sufficient
condition of optimality can be formulated as following.
\begin{thm}
If the necessary condition for the minimum of $f(\cdot)$ at $x_0$
is true and there exists $\beta(g)
>0$ for all suspicious directions $g \in G$ that the inequality
$$
(Ag,g) \meq \beta(g) \| g  \|^2 \,\,\,\,\,\, \forall A \in
\p^2\psi_D f(x_0), \,\,\, \forall D(\cdot) \in \Im,
$$
holds, then $x_0$ is the minimum of $f(\cdot)$.
\label{thmfirstsecsub5}
\end{thm}

To write the necessary and sufficient conditions for the maximum
of $f(\cdot)$ at some point it is sufficiently to remark, that all
maximum points of $f(\cdot)$ are the minimum points of $-
f(\cdot)$.

\section{\bf  Calculus for the subdifferentials of the first and second orders}

Let $f_1,f_2: \mathbb{R}^n \rar \mathbb{R}$ be  Lipschitz
functions and $f(\cdot)=f_1(\cdot)+f_2(\cdot)$. Find for
$f(\cdot), f_1(\cdot), f_2(\cdot)$  the functions $\varphi(\cdot),
\varphi_1(\cdot), \varphi_2(\cdot)$ according to the formulas
written above.

\begin{thm}
The equality \be \p \varphi_D(x_0)=\p \varphi_{1D}(x_0)+ \p
\varphi_{2D}(x_0) \label{firstsecsub22} \ee and the inclusion
$$
\Phi f(x_0) \subset   \Phi f_1(x_0)+\Phi f_2(x_0)
$$
are correct for $f(\cdot)=f_1(\cdot)+f_2(\cdot)$ and any SVM
$D(\cdot) \in \Xi$.
\end{thm}
{\bf Proof.} The equality $f'(z)=f'_1(z)+f'_2(z)$ is true at any
point $z$  of differentiability. As a result, we have
$$
\frac{1}{\mu(D(x))} \int_{D(x)} f'(x+y)d y=\frac{1}{\mu(D(x))}
\int_{D(x)} f'_1(x+y)d y + \frac{1}{\mu(D(x))} \int_{D(x)}
f'_2(x+y)d y.
$$
We get from here (\ref{firstsecsub22}).  Then
$$
\Phi f(x_0)=\bco \, \bigcup_{D(\cdot)} \, \p \varphi_D (x_0)
\subset \bco \bigcup_{D(\cdot)} \, \p \varphi_{1D} (x_0) + \bco
\bigcup_{D(\cdot)} \, \p \varphi_{2D} (x_0) \subset \Phi f_1(x_0)
+ \Phi f_2(x_0).
$$
The theorem is proved. $\Box$

Let  $f(\cdot)=f_1(\cdot)f_2(\cdot)$ now.
\begin{thm}
The equality \be \p \varphi_D(x_0)=\p \varphi_{1D}(x_0)f_2(x_0)+\p
\varphi_{2D}(x_0)f_1(x_0) \label{firstsecsub24} \ee and the
inclusion
$$
\Phi f(x_0) \subset \Phi f_1(x_0) f_2(x_0)+\Phi f_2(x_0) f_1(x_0)
$$
are correct for $f(\cdot)=f_1(\cdot)f_2(\cdot)$ and any SVM
$D(\cdot) \in \Xi$.
\end{thm}
{\bf Proof.} The equality $f'(z)=f'_1(z)f_2(z)+f'_2(z)f_1(z)$ is
true for any point $z$ of differentiability. As a result for SVM
$D(\cdot) \in \Xi$ the equalities
$$
\frac{1}{\mu(D(x))} \int_{D(x)} f'(x+y)d y=\frac{1}{\mu(D(x))}
\int_{D(x)} f'_1(x+y) f_2(x+y) d y +
$$
$$
+\frac{1}{\mu(D(x))} \int_{D(x)} f'_2(x+y) f_1(x+y)d y=
\frac{1}{\mu(D(x))} \int_{D(x)} f'_1(x+y) f_2(x_0) d y +
$$
$$
+\frac{1}{\mu(D(x))} \int_{D(x)} f'_1(x+y) (f_2(x+y)-f_2(x_0)) d y
+ \frac{1}{\mu(D(x))} \int_{D(x)} f'_2(x+y) f_1(x_0) d y +
$$
\be
 +\frac{1}{\mu(D(x))} \int_{D(x)} f'_2(x+y) (f_1(x+y)-f_1(x_0))
dy. \label{firstsecsub23} \ee are correct. From continuity of
$f_1, f_2$ and boundedness of $f'_1, f'_2$  we have for $x \rar
x_0$:
$$
\p \varphi_D(x_0)=\lim_{x \rar x_0} \left[\frac{1}{\mu(D(x))}
\int_{D(x)} f'_1(x+y)d y\right] f_2(x_0) + \lim_{x \rar x_0}
\left[\frac{1}{\mu(D(x))} \int_{D(x)} f'_2(x+y)d y\right] f_1(x_0)
=
$$
\be =\p \varphi_{1D}(x_0)f_2(x_0)+\p \varphi_{2D}(x_0)f_1(x_0).
\label{firstsecsub24} \ee

As soon as (\ref{firstsecsub24}) is true for any SVM $D(\cdot)$,
taking the unity over all $D(\cdot) \in \Xi$ in both sides of the
equality (\ref{firstsecsub24}), we get the statement of the
theorem. $\Box$
\begin{cor}
The equality
$$
\Phi f(x_0)= k \Phi f_1(x_0)
$$
is correct for the Lipschitz function $f(\cdot)=k f_1(\cdot)$,
where $k$ is any constant.
\end{cor}
Pass to the calculus of the subdifferentials of the second order.
Let $f(\cdot)=f_1(\cdot)+f_2(\cdot)$, for which we construct the
functions $\varphi(\cdot), \varphi_1(\cdot), \varphi_2(\cdot)$. We
have for any point $z$, where the matrices
$\varphi''(\cdot),\varphi_1''(\cdot), \varphi_2''(\cdot)$ exist,
and any SVM $D(\cdot) \in \Xi $
$$
\frac{1}{\mu(D(x))} \int_{D(x)} \varphi''(x+y)d
y=\frac{1}{\mu(D(x))} \int_{D(x)} \varphi''_1(x+y)d y +
\frac{1}{\mu(D(x))} \int_{D(x)} \varphi''_2(x+y)d y.
$$
Going to the limit as $x \rar x_0$, we get \be \p^2 \psi_D(x_0)=
\p^2 \psi_{1D}(x_0) + \p^2 \psi_{2D}(x_0), \label{firstsecsub25}
\ee where $\p^2 \psi_{1D}(\cdot)$, $ \p^2 \psi_{2D}(\cdot)$ are
constructed for $\varphi_1(\cdot), \varphi_2(\cdot)$
correspondingly. As soon as (\ref{firstsecsub25}) is correct for
any SVM $D(\cdot)$, then, taking unity over all $D(\cdot) \in \Xi$
in both sides of \ref{firstsecsub25}), we get the statement of the
following theorem.
\begin{thm}
The equality (\ref{firstsecsub25}) and the inclusion
$$
\Psi^2 f(x_0) \incl \Psi^2 f_1(x_0)+\Psi^2 f_2(x_0)
$$
are correct for $f(\cdot)=f_1(\cdot)+ f_2(\cdot)$ and any point
$x_0$.
\end{thm}
Consider the case when $f(\cdot)=f_1(\cdot)f_2(\cdot)$. For
comparison for twice differentiable functions we have
$$
f''(x)=f_1''(x)f_2(x)+f_2''(x)f_1(x)+(f'_1(x))^T f_2'(x) +
(f'_2(x))^T f_1'(x).
$$
Here $(f'_1(x))^T, (f'_2(x))^T $ are the column-vectors received
from the row-vectors $f'_1(x), f'_2(x)$ correspondingly. Go to the
general case. Suppose, that the sets $\Psi^2 f_{1}(x_0)$ and
$\Psi^2 f_{2}(x_0)$ are bounded.

Let us differentiate with respect to $x$ and then take the Steklov
integral from both sides of (\ref{firstsecsub23}) for any SVM
$D(\cdot) \in \Xi$. We get in the result as $x \rar x_0$
$$
\lim_{x \rar x_0} \frac{1}{\mu(D(x))} \int_{D(x)} \varphi''(x+y)d
y= \lim_{x \rar x_0} \left[ \frac{1}{\mu(D(x))} \int_{D(x)}
\varphi''_1(x+y)d y \right] f_2(x_0) +
$$
$$
+ \lim_{x \rar x_0} \left[ \frac{1}{\mu(D(x))} \int_{D(x)}
\varphi''_2(x+y)d y \right] f_1(x_0) + \lim_{x \rar x_0} \left[
\frac{1}{\mu(D(x))} \int_{D(x)} (f'_1(x+y))^T f'_2(x+y) d y
\right] +
$$
$$
+\lim_{x \rar x_0} \left[\frac{1}{\mu(D(x))} \int_{D(x)}
(f'_2(x+y))^T f'_1(x+y) d y \right] + \mbox{terms going to zero as
$x \rar x_0$}.
$$
Really, the missing terms are not  bigger
$$
\lim_{x \rar x_0}  \left[ \| \frac{1}{\mu(D(x))} \int_{D(x)}
\varphi''_1(x+y)  d y  \| \right] \varepsilon_{1}(x) ,
$$
and
$$
\lim_{x \rar x_0} \left[ \| \frac{1}{\mu(D(x))} \int_{D(x)}
\varphi''_2(x+y)  d y \| \right] \varepsilon_{2}(x),
$$
where
$$
\varepsilon_{k}(x)=\max_{y \in D(x)} \vl f_k(x+y)-f_k(x_0) \vl,
\,\, k=1,2.
$$
According to the supposition, the limits of the square brackets
are bounded. Consequently, the limits of whole expressions as $x
\rar x_0$ and the missing terms are equal to zero. We have from
the written expressions  \be \p^2 \psi_D(x_0)=(\p^2
\psi_{1D}(x_0)) f_2(x_0) + (\p^2 \psi_{2D}(x_0)) f_2(x_0) +
{\psi}^2_{12,D}(x_0) + {\psi}^2_{21,D}(x_0), \label{firstsecsub26}
\ee  where
$$
{\psi}^2_{12,D}(x_0) = \{ A \in \mathbb{R}^{n \times n} \vl
\exists\{ x_i \}, A= \lim_{x_i \rar x_0}  \left[
\frac{1}{\mu(D(x_i))} \int_{D(x_i)} (f'_1(x_i+y))^T f'_2(x_i+y) d
y \right],
$$
$$
{\psi}^2_{21,D}(x_0) = \{ A \in \mathbb{R}^{n \times n} \vl
\exists\{ x_i \}, A= \lim_{x_i \rar x_0}  \left[
\frac{1}{\mu(D(x_i))} \int_{D(x_i)} (f'_2(x_i+y))^T f'_1(x_i+y) d
y \right],
$$
the points $x_i$ are taken from the regions of constancy of SVM
$D(\cdot)$.  As a result, we get the following theorem.
\begin{thm}
Under condition  of boundedness of the sets $\Psi^2 f_{1}(x_0)$
and $\Psi^2 f_{2}(x_0)$ the equality (\ref{firstsecsub26}) and the
inclusion
$$
\Psi^2 f(x_0) \incl (\Psi^2 f_1(x_0)) f_2(x_0) + (\Psi^2 f_2(x_0))
f_1(x_0) + {\Psi}^2_{12}(x_0) + {\Psi}^2_{21}(x_0)
$$
hold for $f(\cdot)=f_1(\cdot)f_2(\cdot)$ and any SVM  $D(\cdot)
\in \Xi$ where
$$
{\Psi}^2_{12}(x_0) = \bigcup_{D(\cdot) \in \Xi}
{\psi}^2_{12,D}(x_0), \,\,\, {\Psi}^2_{21}(x_0) =
\bigcup_{D(\cdot) \in \Xi} {\psi}^2_{21,D}(x_0).
$$
\end{thm}

\begin{ex}
Let be $f(x) = \vl x \vl,  x \in \mathbb{R}$. Then $\Phi f(0)=
Df(0)= [-1, 1]=\p_{CL} f(0)$.  The functions $\psi(\cdot)$ are
convex for any constant SVM $D(\cdot)$ according to the qualities
of such functions proved before. During decreasing of the
diameters of the images $D(x)$ the functions $\psi(\cdot),
\psi'(\cdot)$  tend to the functions $f(\cdot), f'(\cdot)$
uniformly on any compact set. Therefore, the second derivatives
$\psi''(x)$ tend to $+\infty$ when $x \rar 0$. From here we have
$\Psi^2 f(0) = \{ +\infty \} .$
\end{ex}

\begin{ex}
Let be $f(\cdot): \mathbb{R} \rar \mathbb{R}$ with a graph lying
between two curves $y=x^2$ and $y=-x^2$ and consisting from
slopes $\pm 1$ with  the limit point  at zero. Then
$$
Df(0)=\{ 0 \}, \,\,\,\, \Psi^2 f(0)= [-2,2].
$$
We can conclude from here that the point zero is not the optimal
point.
\end{ex}

\hspace{0.5cm}

{ \bf Summary}

\hspace{0.2cm}

Proudnikov I.M.

\hspace{0.2cm}

{\large \bf THE SUBDIFFERENTIALS OF THE FIRST AND SECOND ORDERS
FOR LIPSCHITZ FUNCTIONS}

\hspace{0.2cm}

The generalized gradients and matrices of Lipschitz functions  are
defined with the help of the Steklov integral. The
subdifferentials of the first and second orders consisting from
the generalized gradients and the matrices of a Lipschitz function
$f(\cdot)$ are introduced. The Steklov integral over the defined
set-valued mappings is used for the constructions. It is proved
that the subdifferential of the first order coincides with the
average limit values of the integrals of the gradients of
$f(\cdot)$ calculated along curves from a set, introduced earlier
by the author in \cite{lowapp2}. It is proved that the
subdifferentials of the first and second orders are equal to the
first and second derivatives of $f(\cdot)$ correspondingly if such
derivatives exist. The subdifferentials of the first and second
orders are used for formulation of the necessary and sufficient
conditions of optimality.

\end{document}